\documentclass{amsart}
\usepackage{amssymb}
\usepackage{amsfonts}

\newtheorem{theorem}{Theorem}
\theoremstyle{plain}

\newtheorem{corollary}{Corollary}

\newtheorem{definition}{Definition}

\newtheorem{lemma}{Lemma}

\newtheorem{remark}{Remark}

\numberwithin{equation}{section}
\input{tcilatex}

\begin{document}
\title[Divergence Measure for Monotonic Functions]{A General Divergence
Measure for Monotonic Functions and Applications in Information Theory}
\author{S.S. Dragomir}
\address{School of Computer Science and Mathematics\\
Victoria University of Technology\\
PO Box 14428, MCMC 8001\\
Victoria, Australia.}
\email{sever.dragomir@vu.edu.au}
\urladdr{http://rgmia.vu.edu.au/SSDragomirWeb.html}
\date{13 January, 2004.}

\begin{abstract}
A general divergence measure for monotonic functions is introduced. Its
connections with the $f-$divergence for convex functions are explored. The
main properties are pointed out.
\end{abstract}

\maketitle

\section{Introduction}

Let $\left( X,\mathcal{A}\right) $ be a measurable space satisfying $%
\left\vert \mathcal{A}\right\vert >2$ and $\mu $ be a $\sigma -$finite
measure on $\left( X,\mathcal{A}\right) .$ Let $\mathcal{P}$ be the set of
all probability measures on $\left( X,\mathcal{A}\right) $ which are
absolutely continuous with respect to $\mu .$ For $P,Q\in \mathcal{P}$, let $%
p=\frac{dP}{d\mu }$ and $q=\frac{dQ}{d\mu }$ denote the \textit{%
Radon-Nikodym }derivatives of $P$ and $Q$ with respect to $\mu .$

Two probability measures $P,Q\in \mathcal{P}$ are said to be \textit{%
orthogonal} and we denote this by $Q\perp P$ if%
\begin{equation*}
P\left( \left\{ q=0\right\} \right) =Q\left( \left\{ p=0\right\} \right) =1.
\end{equation*}

Let $f:[0,\infty )\rightarrow (-\infty ,\infty ]$ be a convex function that
is continuous at $0,$ i.e., $f\left( 0\right) =\lim_{u\downarrow 0}f\left(
u\right) .$

In 1963, I. Csisz\'{a}r \cite{IC1} introduced the concept of $f-$divergence
as follows.

\begin{definition}
\label{d1.1}Let $P,Q\in \mathcal{P}$. Then%
\begin{equation}
I_{f}\left( Q,P\right) =\int_{X}p\left( x\right) f\left[ \frac{q\left(
x\right) }{p\left( x\right) }\right] d\mu \left( x\right) ,  \label{1.1}
\end{equation}%
is called the $f-$divergence of the probability distributions $Q$ and $P.$
\end{definition}

We now give some examples of $f-$divergences that are well-known and often
used in the literature (see also \cite{CDO}).

\subsection{The Class of $\protect\chi ^{\protect\alpha }-$Divergences}

The $f-$divergences of this class, which is generated by the function $\chi
^{\alpha },$ $\alpha \in \lbrack 1,\infty ),$ defined by%
\begin{equation*}
\chi ^{\alpha }\left( u\right) =\left\vert u-1\right\vert ^{\alpha },\ \ \
u\in \lbrack 0,\infty )
\end{equation*}%
have the form%
\begin{equation}
I_{f}\left( Q,P\right) =\int_{X}p\left\vert \frac{q}{p}-1\right\vert
^{\alpha }d\mu =\int_{X}p^{1-\alpha }\left\vert q-p\right\vert ^{\alpha
}d\mu .  \label{1.2}
\end{equation}%
From this class only the parameter $\alpha =1$ provides a distance in the
topological sense, namely the \textit{total variation distance }$V\left(
Q,P\right) =\int_{X}\left\vert q-p\right\vert d\mu .$ The most prominent
special case of this class is, however, \textit{Karl Pearson's }$\chi ^{2}-$%
\textit{divergence.}

\subsection{Dichotomy Class}

From this class, generated by the function $f_{\alpha }:[0,\infty
)\rightarrow \mathbb{R}$%
\begin{equation*}
f_{\alpha }\left( u\right) =\left\{ 
\begin{array}{ll}
u-1-\ln u & \text{for\ \ }\alpha =0; \\ 
&  \\ 
\frac{1}{\alpha \left( 1-\alpha \right) }\left[ \alpha u+1-\alpha -u^{\alpha
}\right] & \text{for\ \ }\alpha \in \mathbb{R}\backslash \left\{ 0,1\right\}
; \\ 
&  \\ 
1-u+u\ln u & \text{for\ \ }\alpha =1;%
\end{array}%
\right.
\end{equation*}%
only the parameter $\alpha =\frac{1}{2}$ $\left( f_{\frac{1}{2}}\left(
u\right) =2\left( \sqrt{u}-1\right) ^{2}\right) $ provides a distance,
namely, the \textit{Hellinger distance}%
\begin{equation*}
H\left( Q,P\right) =\left[ \int_{X}\left( \sqrt{q}-\sqrt{p}\right) ^{2}d\mu %
\right] ^{\frac{1}{2}}.
\end{equation*}

Another important divergence is the \textit{Kullback-Leibler divergence}
obtained for $\alpha =1,$%
\begin{equation*}
KL\left( Q,P\right) =\int_{X}q\ln \left( \frac{q}{p}\right) d\mu .
\end{equation*}

\subsection{Matsushita's Divergences}

The elements of this class, which is generated by the function $\varphi
_{\alpha },$ $\alpha \in (0,1]$ given by%
\begin{equation*}
\varphi _{\alpha }\left( u\right) :=\left\vert 1-u^{\alpha }\right\vert ^{%
\frac{1}{\alpha }},\ \ \ u\in \lbrack 0,\infty ),
\end{equation*}%
are prototypes of metric divergences, providing the distances $\left[
I_{\varphi _{\alpha }}\left( Q,P\right) \right] ^{\alpha }.$

\subsection{Puri-Vineze Divergences}

This class is generated by the functions $\Phi _{\alpha },$ $\alpha \in
\lbrack 1,\infty )$ given by%
\begin{equation*}
\Phi _{\alpha }\left( u\right) :=\frac{\left\vert 1-u\right\vert ^{\alpha }}{%
\left( u+1\right) ^{\alpha -1}},\ \ \ u\in \lbrack 0,\infty ).
\end{equation*}%
It has been shown in \cite{KOV} that, this class provides the distances $%
\left[ I_{\Phi _{\alpha }}\left( Q,P\right) \right] ^{\frac{1}{\alpha }}.$

\subsection{Divergences of Arimoto-type}

This class is generated by the functions%
\begin{equation*}
\Psi _{\alpha }\left( u\right) :=\left\{ 
\begin{array}{ll}
\frac{\alpha }{\alpha -1}\left[ \left( 1+u^{\alpha }\right) ^{\frac{1}{%
\alpha }}-2^{\frac{1}{\alpha }-1}\left( 1+u\right) \right] & \text{for\ \ }%
\alpha \in \left( 0,\infty \right) \backslash \left\{ 1\right\} ; \\ 
&  \\ 
\left( 1+u\right) \ln 2+u\ln u-\left( 1+u\right) \ln \left( 1+u\right) & 
\text{for\ \ }\alpha =1; \\ 
&  \\ 
\frac{1}{2}\left\vert 1-u\right\vert & \text{for\ \ }\alpha =\infty .%
\end{array}%
\right.
\end{equation*}%
It has been shown in \cite{OV} that, this class provides the distances $%
\left[ I_{\Psi _{\alpha }}\left( Q,P\right) \right] ^{\min \left( \alpha ,%
\frac{1}{\alpha }\right) }$ for $\alpha \in \left( 0,\infty \right) $ and $%
\frac{1}{2}V\left( Q,P\right) $ for $\alpha =\infty .$

\section{Some Classes of Normalised Functions}

We denote by $\mathcal{M}^{=\hspace{-1.63mm}\uparrow }\left( [0,\infty
)\right) $ the class of \textit{monotonic\ nondecreasing functions }defined
on $[0,\infty )$ and by $\mathcal{M}s\left( [0,\infty )\right) $ the class
of \textit{measurable functions }on $[0,\infty ).$ We also consider $%
\mathcal{L}e_{1}\left( [0,\infty )\right) $ the class of measurable
functions $g:[0,\infty )\rightarrow \mathbb{R}$ with the property that%
\begin{equation}
g\left( t\right) \leq g\left( 1\right) \leq g\left( s\right) \text{ \ for \ }%
0\leq t\leq 1\leq s<\infty .  \label{2.1}
\end{equation}%
It is obvious that%
\begin{equation}
\mathcal{M}^{=\hspace{-1.63mm}\uparrow }\left( [0,\infty )\right)
\subsetneqq \mathcal{L}e_{1}\left( [0,\infty )\right) ,  \label{2.2}
\end{equation}%
and the inclusion (\ref{2.2}) is strict.

We say that a function $f:[0,\infty )\rightarrow \mathbb{R}$ is \textit{%
normalised }if $f\left( 1\right) =0.$ We denote by $\mathcal{M}s_{0}\left(
[0,\infty )\right) $ the class of all normalised measurable functions
defined on $[0,\infty ).$ We also need the following classes of functions%
\begin{align*}
\mathcal{C}o\left( [0,\infty )\right) & :=\left\{ f\in \mathcal{M}%
s_{0}\left( [0,\infty )\right) |f\text{ is continuous convex on }[0,\infty
)\right\} ; \\
\mathcal{D}_{0}\left( [0,\infty )\right) & :=\left\{ f\in \mathcal{M}%
s_{0}\left( [0,\infty )\right) |f\left( t\right) =\left( t-1\right) g\left(
t\right) ,\ \forall t\in \lbrack 0,\infty ),\ g\in \mathcal{M}^{=\hspace{%
-1.63mm}\uparrow }\left( [0,\infty )\right) \right\} ;
\end{align*}%
and%
\begin{equation*}
\mathcal{O}_{0}\left( [0,\infty )\right) :=\left\{ f\in \mathcal{M}%
s_{0}\left( [0,\infty )\right) |f\left( t\right) =\left( t-1\right) g\left(
t\right) ,\ \forall t\in \lbrack 0,\infty ),\ g\in \mathcal{L}e_{1}\left(
[0,\infty )\right) \right\} .
\end{equation*}%
From the definition of $\mathcal{D}_{0}\left( [0,\infty )\right) $ and $%
\mathcal{O}_{0}\left( [0,\infty )\right) $ and taking into account that the
strict inclusion (\ref{2.2}) holds, we deduce that%
\begin{equation}
\mathcal{D}_{0}\left( [0,\infty )\right) \subsetneqq \mathcal{O}_{0}\left(
[0,\infty )\right) ,  \label{2.3}
\end{equation}%
and the inclusion is strict.

For the other two classes, we may state the following result.

\begin{lemma}
\label{l2.1}We have the strict inclusion%
\begin{equation}
\mathcal{C}o\left( [0,\infty )\right) \subsetneqq \mathcal{D}_{0}\left(
[0,\infty )\right) .  \label{2.4}
\end{equation}
\end{lemma}

\begin{proof}
We will show that any continuous convex function $f:[0,\infty )\rightarrow 
\mathbb{R}$ that is normalised may be represented as:%
\begin{equation}
f\left( t\right) =\left( t-1\right) g\left( t\right) \ \text{for any \ }t\in
\lbrack 0,\infty ),  \label{2.5}
\end{equation}%
where $g\in \mathcal{M}^{=\hspace{-1.63mm}\uparrow }\left( [0,\infty
)\right) .$

Now, let $f\in \mathcal{C}o\left( [0,\infty )\right) .$ For $\lambda \in %
\left[ D_{-}f\left( 1\right) ,D_{+}f\left( 1\right) \right] ,$ define%
\begin{equation*}
g_{\lambda }\left( t\right) :=\left\{ 
\begin{array}{ll}
\dfrac{f\left( t\right) }{t-1} & \text{if\ \ }t\in \lbrack 0,1)\cup \left(
1,\infty \right) , \\ 
&  \\ 
\lambda & \text{if\ \ }t=1.%
\end{array}%
\right.
\end{equation*}%
We use the following well known result \cite[p. 111]{GB}:

If $\Psi $ is convex on $\left( a,b\right) $ and $a<s<t<u<b,$ then%
\begin{equation}
\Psi \left( s,t\right) \leq \Psi \left( s,u\right) \leq \Psi \left(
t,u\right) ,  \label{2.6}
\end{equation}%
where%
\begin{equation*}
\Psi \left( s,t\right) =\frac{\Psi \left( t\right) -\Psi \left( s\right) }{%
t-s}.
\end{equation*}%
If $\Psi $ is strictly convex on $\left( a,b\right) ,$ equality will not
occur in (\ref{2.6}).

If we apply the above result for $0<s<t<1,$ then we can state%
\begin{equation*}
\dfrac{f\left( s\right) }{s-1}\leq \dfrac{f\left( t\right) }{t-1}.
\end{equation*}%
Taking the limit over $t\rightarrow 1,$ $t<1,$ we deduce%
\begin{equation*}
\dfrac{f\left( s\right) }{s-1}\leq D_{-}f\left( 1\right)
\end{equation*}%
showing that for $0<t<1,$ we have $g_{\lambda }\left( t\right) \leq \lambda
. $

Similarly, we may prove that for $1<t<\infty ,$ $g_{\lambda }\left( t\right)
\geq \lambda .$ If we use the same result for $0<t_{1}<t_{2}<1,$ then we may
write%
\begin{equation*}
\dfrac{f\left( t_{1}\right) }{t_{1}-1}\leq \dfrac{f\left( t_{2}\right) }{%
t_{2}-1},
\end{equation*}%
which gives $g_{\lambda }\left( t_{1}\right) \leq g_{\lambda }\left(
t_{2}\right) $ for $0<t_{1}<t_{2}<1.$

In a similar fashion we can prove that for $1<t_{1}<t_{2}<\infty ,$ $%
g_{\lambda }\left( t_{1}\right) \leq g_{\lambda }\left( t_{2}\right) ,$ and
thus we may conclude that the function $g_{\lambda }$ is monotonic
non-decreasing on the whole interval $[0,\infty ).$

If we consider now the function $f\left( t\right) =\left( t-1\right) e^{\eta
t},$ $t\in \lbrack 0,\infty ),$ we observe that $f^{\prime }\left( t\right)
=\left( \eta t-3\right) e^{\eta t},$ $f^{\prime \prime }\left( t\right)
=8e^{\eta t}\left( 2t-1\right) $ which shows that $f$ is not convex on $%
[0,\infty ).$ Obviously, $f\in \mathcal{D}_{0}\left( [0,\infty )\right) ,$
and thus the inclusion (\ref{2.4}) is indeed strict.
\end{proof}

\begin{remark}
\label{r2.1}If $f\in \mathcal{D}_{0}\left( [0,\infty )\right) $ and $%
g_{1},g_{2}\in \mathcal{M}^{=\hspace{-1.63mm}\uparrow }\left( [0,\infty
)\right) $ are two functions with%
\begin{equation*}
f\left( t\right) =\left( t-1\right) g_{1}\left( t\right) ,\ \ \ f\left(
t\right) =\left( t-1\right) g_{2}\left( t\right)
\end{equation*}%
for each $t\in \lbrack 0,\infty ),$ then we get%
\begin{equation*}
\left( t-1\right) \left[ g_{1}\left( t\right) -g_{2}\left( t\right) \right]
=0
\end{equation*}%
for any $t\in \lbrack 0,\infty )$ showing that $g_{1}\left( t\right)
=g_{2}\left( t\right) $ for each $t\in \lbrack 0,1)\cup \left( 1,\infty
\right) .$ They may have different values in $t=1$.
\end{remark}

\section{Some Fundamental Properties of $f-$Divergence for $f\in \mathcal{C}%
o\left( [0,\infty )\right) $}

For $f\in \mathcal{C}o\left( [0,\infty )\right) $ we obtain the $\ast -$%
\textit{conjugate }function of $f$ by%
\begin{equation*}
f^{\ast }\left( u\right) =uf\left( \frac{1}{u}\right) ,\ \ \ u\in \left(
0,\infty \right) .
\end{equation*}%
It is also known that if $f\in \mathcal{C}o\left( [0,\infty )\right) ,$ then 
$f^{\ast }\in \mathcal{C}o\left( [0,\infty )\right) .$

The following two theorems contain the most basic properties of $f-$%
divergences. For their proof we refer the reader to Chapter 1 of \cite{LV}
(see also \cite{CDO}).

\begin{theorem}[Uniqueness and Symmetry Theorem]
\label{t3.1}Let $f,f_{1}$ be continuous convex on $[0,\infty ).$

\begin{enumerate}
\item[(i)] We have%
\begin{equation*}
I_{f_{1}}\left( Q,P\right) =I_{f}\left( Q,P\right) ,
\end{equation*}%
for any $P,Q\in \mathcal{P}$ if and only if there exists a constant $c\in 
\mathbb{R}$ such that%
\begin{equation*}
f_{1}\left( u\right) =f\left( u\right) +c\left( u-1\right) ,
\end{equation*}%
for any $u\in \lbrack 0,\infty );$

\item[(ii)] We have%
\begin{equation*}
I_{f^{\ast }}\left( Q,P\right) =I_{f}\left( Q,P\right) ,
\end{equation*}%
for any $P,Q\in \mathcal{P}$ if and only if there exists a constant $d\in 
\mathbb{R}$ such that%
\begin{equation*}
f^{\ast }\left( u\right) =f\left( u\right) +d\left( c-1\right) ,
\end{equation*}%
for any $u\in \lbrack 0,\infty ).$
\end{enumerate}
\end{theorem}

\begin{theorem}[Range of Values Theorem]
\label{t3.2}Let $f:[0,\infty )\rightarrow \mathbb{R}$ be a continuous convex
function on $[0,\infty ).$

For any $P,Q\in \mathcal{P}$, we have the double inequality%
\begin{equation}
f\left( 1\right) \leq I_{f}\left( Q,P\right) \leq f\left( 0\right) +f^{\ast
}\left( 0\right) .  \label{3.1}
\end{equation}

\begin{enumerate}
\item[(i)] If $P=Q,$ then the equality holds in the first part of (\ref{3.1}%
).

If $f$ is strictly convex at $1,$ then the equality holds in the first part
of (\ref{3.1}) if and only if $P=Q;$

\item[(ii)] If $Q\perp P,$ then the equality holds in the second part of (%
\ref{3.1}).

If $f\left( 0\right) +f^{\ast }\left( 0\right) <\infty ,$ then equality
holds in the second part of (\ref{3.1}) if and only if $Q\perp P.$
\end{enumerate}
\end{theorem}

Define the function $\tilde{f}:\left( 0,\infty \right) \rightarrow \mathbb{R}
$, $\tilde{f}\left( u\right) =\frac{1}{2}\left( f\left( u\right) +f^{\ast
}\left( u\right) \right) .$ The following result is a refinement of the
second inequality in Theorem \ref{t3.2} (see \cite[Theorem 3]{CDO}).

\begin{theorem}
\label{t3.3}Let $f\in \mathcal{C}o\left( [0,\infty )\right) $ with $f\left(
0\right) +f^{\ast }\left( 0\right) <\infty .$ Then%
\begin{equation}
0\leq I_{f}\left( Q,P\right) \leq \tilde{f}\left( 0\right) V\left( Q,P\right)
\label{3.2}
\end{equation}%
for any $Q,P\in \mathcal{P}$.
\end{theorem}

\section{A General Divergence Measure}

If $f:[0,\infty )\rightarrow \mathbb{R}$ is a general measurable function,
then we may define the $f-$\textit{divergence} in the same way, i.e., if $%
P,Q\in \mathcal{P}$, then%
\begin{equation*}
I_{f}\left( Q,P\right) =\int_{X}p\left( x\right) f\left[ \frac{q\left(
x\right) }{p\left( x\right) }\right] d\mu \left( x\right) .
\end{equation*}%
For a measurable function $g:[0,\infty )\rightarrow \mathbb{R}$, we may also
define the $\delta -$\textit{divergence} by the formula%
\begin{equation*}
\delta _{g}\left( Q,P\right) =\int_{X}\left[ q\left( x\right) -p\left(
x\right) \right] g\left[ \frac{q\left( x\right) }{p\left( x\right) }\right]
d\mu \left( x\right) .
\end{equation*}

It is obvious that the $\delta -$divergence of a function $g$ may be seen as
the $f-$divergence of the function $f,$ where $f\left( t\right) =\left(
t-1\right) g\left( t\right) $ for $t\in \lbrack 0,\infty ).$

If $f\in \mathcal{C}o\left( [0,\infty )\right) $ and since $f\left( t\right)
=\left( t-1\right) g_{\lambda }\left( t\right) ,$ $t\in \lbrack 0,\infty ),$
we have%
\begin{equation}
g_{\lambda }\left( t\right) :=\left\{ 
\begin{array}{ll}
\dfrac{f\left( t\right) }{t-1} & \text{if\ \ }t\in \lbrack 0,1)\cup \left(
1,\infty \right) , \\ 
&  \\ 
\lambda  & \text{if\ \ }t=1;%
\end{array}%
\right.   \label{4.0}
\end{equation}%
and $\lambda \in \left[ D_{-}f\left( 1\right) ,D_{+}f\left( 1\right) \right]
,$ shows that for any $f\in \mathcal{C}o\left( [0,\infty )\right) $ we have%
\begin{equation}
I_{f}\left( Q,P\right) =\delta g_{\lambda }\left( Q,P\right) \ \ \text{for
any \ }P,Q\in \mathcal{P},  \label{4.1}
\end{equation}%
i.e.,\textit{\ the }$f-$\textit{divergence for any normalised continuous
convex function }$f:[0,\infty )\rightarrow \mathbb{R}$\textit{\ may be seen
as the }$\delta -$\textit{divergence of the function }$g_{\lambda }$\textit{%
\ defined by} (\ref{4.0}).

In what follows, we point out some fundamental properties of the $\delta -$%
divergence.

\begin{theorem}
\label{t4.1}Let $g:[0,\infty )\rightarrow \mathbb{R}$ be a measurable
function on $[0,\infty )$ and $P,Q\in \mathcal{P}$. If there exists the
constants $m,M$ with%
\begin{equation}
-\infty <m\leq g\left[ \frac{q\left( x\right) }{p\left( x\right) }\right]
\leq M<\infty  \label{4.2}
\end{equation}%
for $\mu -$a.e. $x\in X,$ then we have the inequality%
\begin{equation}
\left\vert \delta _{g}\left( Q,P\right) \right\vert \leq \frac{1}{2}\left(
M-m\right) V\left( Q,P\right) .  \label{4.3}
\end{equation}
\end{theorem}

\begin{proof}
We observe that the following identity holds true%
\begin{equation}
\delta _{g}\left( Q,P\right) =\int_{X}\left[ q\left( x\right) -p\left(
x\right) \right] \left[ g\left[ \frac{q\left( x\right) }{p\left( x\right) }%
\right] -\frac{m+M}{2}\right] d\mu \left( x\right)  \label{4.4}
\end{equation}%
By (\ref{4.2}), we deduce that%
\begin{equation*}
\left\vert g\left[ \frac{q\left( x\right) }{p\left( x\right) }\right] -\frac{%
m+M}{2}\right\vert \leq \frac{1}{2}\left( M-m\right)
\end{equation*}%
for $\mu -$a.e. $x\in X.$

Taking the modulus in (\ref{4.4}) we deduce%
\begin{align*}
\left\vert \delta _{g}\left( Q,P\right) \right\vert & \leq
\int_{X}\left\vert q\left( x\right) -p\left( x\right) \right\vert \left\vert
g\left[ \frac{q\left( x\right) }{p\left( x\right) }-\frac{m+M}{2}\right]
\right\vert d\mu \left( x\right) \\
& \leq \frac{1}{2}\left( M-m\right) \int_{X}\left\vert q\left( x\right)
-p\left( x\right) \right\vert d\mu \left( x\right) \\
& =\frac{1}{2}\left( M-m\right) V\left( Q,P\right)
\end{align*}%
and the inequality (\ref{4.3}) is proved.
\end{proof}

The following corollary is a natural consequence of the above theorem.

\begin{corollary}
\label{c4.2}Let $g:[0,\infty )\rightarrow \mathbb{R}$ be a measurable
function on $[0,\infty )$. If%
\begin{equation*}
m:=ess\inf_{t\in \lbrack 0,\infty )}g\left( t\right) >-\infty ,\ \ \ \
M:=ess\sup_{t\in \lbrack 0,\infty )}g\left( t\right) <\infty ,
\end{equation*}%
then for any $P,Q\in \mathcal{P}$, we have the inequality%
\begin{equation}
\left\vert \delta _{g}\left( Q,P\right) \right\vert \leq \frac{1}{2}\left(
M-m\right) V\left( Q,P\right) .  \label{4.4.a}
\end{equation}
\end{corollary}

\begin{remark}
\label{r4.3}We know that, if $f:[0,\infty )\rightarrow \mathbb{R}$ is a
normalised continuous convex function and if $\lim_{t\downarrow 0}f^{\ast
}\left( t\right) =\lim_{u\downarrow 0}\left[ uf\left( \frac{1}{u}\right) %
\right] =:f^{\ast }\left( 0\right) ,$ then we have the inequality [Theorem
2.3]%
\begin{equation}
I_{f}\left( Q,P\right) \leq \frac{f\left( 0\right) +f^{\ast }\left( 0\right) 
}{2}V\left( Q,P\right) ,  \label{4.5}
\end{equation}%
for any $P,Q\in \mathcal{P}$. We can prove this inequality by the use of
Corollary \ref{c4.2} as follows. We have%
\begin{equation*}
I_{f}\left( Q,P\right) =\delta g_{\lambda }\left( Q,P\right) ,
\end{equation*}%
where%
\begin{equation*}
g_{\lambda }\left( t\right) :=\left\{ 
\begin{array}{ll}
\dfrac{f\left( t\right) }{t-1} & \text{if\ \ }t\in \lbrack 0,1)\cup \left(
1,\infty \right) , \\ 
&  \\ 
\lambda  & \text{if\ \ }t=1,%
\end{array}%
\right. 
\end{equation*}%
where $\lambda \in \left[ D_{-}f\left( 1\right) ,D_{+}f\left( 1\right) %
\right] $ and $g_{\lambda }\in \mathcal{M}^{=\hspace{-1.63mm}\uparrow
}\left( [0,\infty )\right) .$ We observe that for any $t\in \lbrack 0,\infty
),$ we have%
\begin{equation*}
g_{\lambda }\left( t\right) \geq \lim_{t\rightarrow 0+}g_{\lambda }\left(
t\right) =-f\left( 0\right) =m>-\infty 
\end{equation*}%
and%
\begin{align*}
g_{\lambda }\left( t\right) & \leq \lim_{t\rightarrow +\infty }g_{\lambda
}\left( t\right) =\lim_{t\rightarrow +\infty }\dfrac{f\left( t\right) }{t-1}%
=\lim_{u\rightarrow 0+}\left[ \frac{f\left( \frac{1}{u}\right) }{\frac{1}{u}%
-1}\right]  \\
& =\lim_{u\rightarrow 0+}\left[ \frac{uf\left( \frac{1}{u}\right) }{1-u}%
\right] =f^{\ast }\left( 0\right) =M<\infty .
\end{align*}%
Applying Corollary \ref{c4.2} for $m=-f\left( 0\right) $ and $M=f^{\ast
}\left( 0\right) ,$ we deduce the desired inequality (\ref{4.5}).
\end{remark}

The following result also holds.

\begin{theorem}
\label{t4.4}Let $g:[0,\infty )\rightarrow \mathbb{R}$ be a measurable
function on $[0,\infty )$ and $P,Q\in \mathcal{P}$. If there exists a
constant $K$ with $K>0$ such that%
\begin{equation}
\left\vert g\left( \frac{q\left( x\right) }{p\left( x\right) }\right)
-g\left( 1\right) \right\vert \leq K\left\vert \frac{q\left( x\right) }{%
p\left( x\right) }-1\right\vert ^{\alpha },  \label{4.6}
\end{equation}%
for $\mu -$a.e. $x\in X,$ where $\alpha \in \left( 0,\infty \right) $ is a
given number, then we have the inequality%
\begin{equation}
\left\vert \delta _{g}\left( Q,P\right) \right\vert \leq KI_{\chi ^{\alpha
+1}}\left( Q,P\right) .  \label{4.7}
\end{equation}
\end{theorem}

\begin{proof}
We observe that the following identity holds true%
\begin{equation}
\delta _{g}\left( Q,P\right) =\int_{X}\left[ q\left( x\right) -p\left(
x\right) \right] \left[ g\left[ \frac{q\left( x\right) }{p\left( x\right) }%
\right] -g\left( 1\right) \right] d\mu \left( x\right) .  \label{4.8}
\end{equation}%
Taking the modulus in (\ref{4.8}) and using the condition (\ref{4.6}), we
have successively%
\begin{align*}
\left\vert \delta _{g}\left( Q,P\right) \right\vert & \leq
\int_{X}\left\vert q\left( x\right) -p\left( x\right) \right\vert \left\vert
g\left[ \frac{q\left( x\right) }{p\left( x\right) }\right] -g\left( 1\right)
\right\vert d\mu \left( x\right) \\
& \leq K\int_{X}\left[ p\left( x\right) \right] ^{-\alpha }\left\vert
q\left( x\right) -p\left( x\right) \right\vert ^{\alpha +1}d\mu \left(
x\right) \\
& \leq KI_{\chi ^{\alpha +1}}\left( Q,P\right)
\end{align*}%
and the inequality (\ref{4.7}) is obtained.
\end{proof}

The following corollary holds.

\begin{corollary}
\label{c4.5}Let $g:[0,\infty )\rightarrow \mathbb{R}$ be a measurable
function on $[0,\infty )$ with the property that there exists a constant $K$
with the property that%
\begin{equation}
\left\vert g\left( t\right) -g\left( 1\right) \right\vert \leq K\left\vert
t-1\right\vert ^{\alpha },  \label{4.9}
\end{equation}%
for a.e. $t\in \lbrack 0,\infty ),$ where $\alpha >0$ is a given number.
Then for any $P,Q\in \mathcal{P}$, we have the inequality%
\begin{equation}
\left\vert \delta _{g}\left( Q,P\right) \right\vert \leq KI_{\chi ^{\alpha
+1}}\left( Q,P\right) .  \label{4.10}
\end{equation}
\end{corollary}

\begin{remark}
\label{r4.6}If the function $g:[0,\infty )\rightarrow \mathbb{R}$ is H\"{o}%
lder continuous with a constant $H>0$ and $\beta \in (0,1],$ i.e., 
\begin{equation*}
\left\vert g\left( t\right) -g\left( s\right) \right\vert \leq H\left\vert
t-s\right\vert ^{\beta },
\end{equation*}%
for any $t,s\in \lbrack 0,\infty ),$ then obviously (\ref{4.5}) holds with $%
K=H$ and $\alpha =\beta .$

If $g:[0,\infty )\rightarrow \mathbb{R}$ is Lipschitzian with the constant $%
L>0,$ i.e.,%
\begin{equation*}
\left\vert g\left( t\right) -g\left( s\right) \right\vert \leq L\left\vert
t-s\right\vert ,
\end{equation*}%
for any $t,s\in \lbrack 0,\infty ),$ then%
\begin{equation}
\left\vert \delta _{g}\left( Q,P\right) \right\vert \leq KI_{\chi
^{2}}\left( Q,P\right) ,  \label{4.11}
\end{equation}%
for any $P,Q\in \mathcal{P}$.

Finally, if $g$ is locally absolutely continuous and the derivative $%
g^{\prime }:[0,\infty )\rightarrow \mathbb{R}$ is essentially bounded, i.e., 
$\left\Vert g^{\prime }\right\Vert _{[0,\infty ),\infty }:=ess\sup_{t\in
\lbrack 0,\infty )}\left\vert g^{\prime }\left( t\right) \right\vert <\infty
,$ then we have the inequality%
\begin{equation}
\left\vert \delta _{g}\left( Q,P\right) \right\vert \leq \left\Vert
g^{\prime }\right\Vert _{[0,\infty ),\infty }I_{\chi ^{2}}\left( Q,P\right) ,
\label{4.12}
\end{equation}%
for any $P,Q\in \mathcal{P}$.
\end{remark}

The following result concerning $f-$divergences for $f$ convex functions
holds.

\begin{theorem}
\label{t4.6}Let $f:\left[ 0,\infty \right] \rightarrow \mathbb{R}$ be a
continuous convex function on $[0,\infty ).$ If $\lambda \in \left[
D_{-}f\left( 1\right) ,D_{+}f\left( 1\right) \right] $ ($\lambda =f^{\prime
}\left( 1\right) $ if $f$ is differentiable at $t=1$), and there exists a
constant $K>0$ and $\alpha >0$ such that%
\begin{equation}
\left\vert f\left( t\right) -\lambda \left( t-1\right) \right\vert \leq
K\left\vert t-1\right\vert ^{\alpha +1},  \label{4.13}
\end{equation}%
for any $t\in \lbrack 0,\infty ),$ then we have the inequality%
\begin{equation}
0\leq I_{f}\left( Q,P\right) \leq KI_{\chi ^{\alpha +1}}\left( Q,P\right) ,
\label{4.14}
\end{equation}%
for any $P,Q\in \mathcal{P}$.
\end{theorem}

\begin{proof}
We have%
\begin{equation*}
I_{f}\left( Q,P\right) =\int_{X}\left[ q\left( x\right) -p\left( x\right) %
\right] g_{\lambda }\left[ \frac{p\left( x\right) }{q\left( x\right) }\right]
d\mu \left( x\right) =\delta g_{\lambda }\left( Q,P\right) ,
\end{equation*}%
where%
\begin{equation*}
g_{\lambda }\left( t\right) :=\left\{ 
\begin{array}{ll}
\dfrac{f\left( t\right) }{t-1} & \text{if\ \ }t\in \lbrack 0,1)\cup \left(
1,\infty \right) , \\ 
&  \\ 
\lambda & \text{if\ \ }t=1,%
\end{array}%
\right.
\end{equation*}%
and $\lambda \in \left[ D_{-}f\left( 1\right) ,D_{+}f\left( 1\right) \right]
.$

Applying Corollary \ref{c4.5} for $g_{\lambda },$ we deduce the desired
result.
\end{proof}

\section{The Positivity of $\protect\delta -$Divergence for $g\in \mathcal{M}%
^{=\hspace{-1.63mm}\uparrow }\left( [0,\infty )\right) $}

The following result holds.

\begin{theorem}
\label{t5.1}If $g\in \mathcal{M}^{=\hspace{-1.63mm}\uparrow }\left(
[0,\infty )\right) ,$ then $\delta _{g}\left( Q,P\right) \geq 0$ for any $%
P,Q\in \mathcal{P}$.
\end{theorem}

\begin{proof}
We use the identity%
\begin{align}
& \delta _{g}\left( Q,P\right)  \label{5.1} \\
& =\int_{X}\left[ q\left( x\right) -p\left( x\right) \right] g\left[ \frac{%
q\left( x\right) }{p\left( x\right) }\right] d\mu \left( x\right)  \notag \\
& =\int_{X}p\left( x\right) \left[ \frac{q\left( x\right) }{p\left( x\right) 
}-1\right] g\left[ \frac{q\left( x\right) }{p\left( x\right) }\right] d\mu
\left( x\right)  \notag \\
& =\frac{1}{2}\int_{X}\int_{X}p\left( x\right) p\left( y\right) \left[ \frac{%
q\left( x\right) }{p\left( x\right) }-\frac{q\left( y\right) }{p\left(
y\right) }\right] \left[ g\left[ \frac{q\left( x\right) }{p\left( x\right) }%
\right] -g\left[ \frac{q\left( y\right) }{p\left( y\right) }\right] \right]
d\mu \left( x\right) d\mu \left( y\right) .  \notag
\end{align}%
Since $g\in \mathcal{M}^{=\hspace{-1.63mm}\uparrow }\left( [0,\infty
)\right) ,$ then for any $t,s\in \lbrack 0,\infty ),$ we have%
\begin{equation*}
\left( t-s\right) \left( g\left( t\right) -g\left( s\right) \right) \geq 0
\end{equation*}%
giving that%
\begin{equation*}
\left[ \frac{q\left( x\right) }{p\left( x\right) }-\frac{q\left( y\right) }{%
p\left( y\right) }\right] \left[ g\left[ \frac{q\left( x\right) }{p\left(
x\right) }\right] -g\left[ \frac{q\left( y\right) }{p\left( y\right) }\right]
\right] \geq 0
\end{equation*}%
for any $x,y\in X.$

Using the representation (\ref{5.1}), we deduce the desired result.
\end{proof}

The following corollary is a natural consequence of the above result.

\begin{corollary}
\label{c5.2}If $f\in \mathcal{D}_{0}\left( [0,\infty )\right) ,$ then $%
I_{f}\left( Q,P\right) \geq 0$ for any $P,Q\in \mathcal{P}$.
\end{corollary}

\begin{proof}
If $f\in \mathcal{D}_{0}\left( [0,\infty )\right) ,$ then there exists a $%
g\in \mathcal{M}^{=\hspace{-1.63mm}\uparrow }\left( [0,\infty )\right) $
such that $f\left( t\right) =\left( t-1\right) g\left( t\right) $ for any $%
t\in \lbrack 0,\infty ).$ Then%
\begin{align*}
I_{f}\left( Q,P\right) & =\int_{X}p\left( x\right) f\left[ \frac{q\left(
x\right) }{p\left( x\right) }\right] d\mu \left( x\right) \\
& =\int_{X}p\left( x\right) \left[ \frac{q\left( x\right) }{p\left( x\right) 
}-1\right] g\left[ \frac{q\left( x\right) }{p\left( x\right) }\right] d\mu
\left( x\right) \\
& =\delta _{g}\left( Q,P\right) \geq 0,
\end{align*}%
and the proof is completed.
\end{proof}

In fact, the following improvement of Theorem \ref{t5.1} holds.

\begin{theorem}
\label{t5.3}If $g\in \mathcal{M}^{=\hspace{-1.63mm}\uparrow }\left(
[0,\infty )\right) ,$ then%
\begin{equation}
\delta _{g}\left( Q,P\right) \geq \left\vert \delta _{\left\vert
g\right\vert }\left( Q,P\right) \right\vert \geq 0,  \label{5.2}
\end{equation}%
for any $P,Q\in \mathcal{P}$.
\end{theorem}

\begin{proof}
Since $g$ is monotonic nondecreasing, we have%
\begin{align}
& \left[ \frac{q\left( x\right) }{p\left( x\right) }-\frac{q\left( y\right) 
}{p\left( y\right) }\right] \left[ g\left[ \frac{q\left( x\right) }{p\left(
x\right) }\right] -g\left[ \frac{q\left( y\right) }{p\left( y\right) }\right]
\right]  \label{5.3} \\
& =\left\vert \left( \frac{q\left( x\right) }{p\left( x\right) }-\frac{%
q\left( y\right) }{p\left( y\right) }\right) \left( g\left[ \frac{q\left(
x\right) }{p\left( x\right) }\right] -g\left[ \frac{q\left( y\right) }{%
p\left( y\right) }\right] \right) \right\vert  \notag \\
& \geq \left\vert \left( \frac{q\left( x\right) }{p\left( x\right) }-\frac{%
q\left( y\right) }{p\left( y\right) }\right) \left( \left\vert g\left[ \frac{%
q\left( x\right) }{p\left( x\right) }\right] \right\vert -\left\vert g\left[ 
\frac{q\left( y\right) }{p\left( y\right) }\right] \right\vert \right)
\right\vert  \notag
\end{align}%
for any $x,y\in X.$

Multiplying (\ref{5.3}) by $p\left( x\right) p\left( y\right) \geq 0$ and
integrating on $X^{2},$ we deduce%
\begin{multline*}
\int_{X}\int_{X}p\left( x\right) p\left( y\right) \left( \frac{q\left(
x\right) }{p\left( x\right) }-\frac{q\left( y\right) }{p\left( y\right) }%
\right) \left[ g\left[ \frac{q\left( x\right) }{p\left( x\right) }\right] -g%
\left[ \frac{q\left( y\right) }{p\left( y\right) }\right] \right] d\mu
\left( x\right) d\mu \left( y\right) \\
\geq \left\vert \int_{X}\int_{X}p\left( x\right) p\left( y\right) \left( 
\frac{q\left( x\right) }{p\left( x\right) }-\frac{q\left( y\right) }{p\left(
y\right) }\right) \left( g\left[ \frac{q\left( x\right) }{p\left( x\right) }%
\right] -g\left[ \frac{q\left( y\right) }{p\left( y\right) }\right] \right)
d\mu \left( x\right) d\mu \left( y\right) \right\vert .
\end{multline*}%
Using the representation (\ref{5.1}) and the same identity for $\left\vert
g\right\vert ,$ we deduce the desired inequality (\ref{5.2}).
\end{proof}

Before we point out other possible refinements for the positivity inequality 
$\delta _{g}\left( Q,P\right) \geq 0,$ where $g\in \mathcal{M}^{=\hspace{%
-1.63mm}\uparrow }\left( [0,\infty )\right) ,$ we need the following
divergence measure as well:%
\begin{equation*}
\bar{\delta}_{h}\left( Q,P\right) :=\int_{X}\left\vert q\left( x\right)
-p\left( x\right) \right\vert h\left[ \frac{q\left( x\right) }{p\left(
x\right) }\right] d\mu \left( x\right)
\end{equation*}%
which will be called the \textit{absolute }$\delta -$\textit{divergence}
generated by the function $h:[0,\infty )\rightarrow \mathbb{R}$ that is
assumed to be measurable on $[0,\infty ).$

The following result holds.

\begin{theorem}
\label{t5.4}If $g\in \mathcal{M}^{=\hspace{-1.63mm}\uparrow }\left(
[0,\infty )\right) ,$ then%
\begin{multline}
\delta _{g}\left( Q,P\right)  \label{5.4} \\
\geq \max \left\{ \left\vert \bar{\delta}_{g}\left( Q,P\right) -V\left(
Q,P\right) I_{g}\left( Q,P\right) \right\vert ,\left\vert \bar{\delta}%
_{\left\vert g\right\vert }\left( Q,P\right) -V\left( Q,P\right)
I_{\left\vert g\right\vert }\left( Q,P\right) \right\vert \right\} \geq 0,
\end{multline}%
for any $P,Q\in \mathcal{P}$.
\end{theorem}

\begin{proof}
Since $g$ is monotonic, we have%
\begin{align}
& \left( \frac{q\left( x\right) }{p\left( x\right) }-\frac{q\left( y\right) 
}{p\left( y\right) }\right) \left( g\left[ \frac{q\left( x\right) }{p\left(
x\right) }\right] -g\left[ \frac{q\left( y\right) }{p\left( y\right) }\right]
\right)  \label{5.5} \\
& =\left\vert \left[ \left( \frac{q\left( x\right) }{p\left( x\right) }%
-1\right) -\left( \frac{q\left( y\right) }{p\left( y\right) }-1\right) %
\right] \left[ g\left[ \frac{q\left( x\right) }{p\left( x\right) }\right] -g%
\left[ \frac{q\left( y\right) }{p\left( y\right) }\right] \right] \right\vert
\notag \\
& \geq \left\{ 
\begin{array}{l}
\left\vert \left[ \left\vert \frac{q\left( x\right) }{p\left( x\right) }%
-1\right\vert -\left\vert \frac{q\left( y\right) }{p\left( y\right) }%
-1\right\vert \right] \left[ g\left[ \frac{q\left( x\right) }{p\left(
x\right) }\right] -g\left[ \frac{q\left( y\right) }{p\left( y\right) }\right]
\right] \right\vert \\ 
\\ 
\left\vert \left[ \left\vert \frac{q\left( x\right) }{p\left( x\right) }%
-1\right\vert -\left\vert \frac{q\left( y\right) }{p\left( y\right) }%
-1\right\vert \right] \left[ \left\vert g\left[ \frac{q\left( x\right) }{%
p\left( x\right) }\right] \right\vert -\left\vert g\left[ \frac{q\left(
y\right) }{p\left( y\right) }\right] \right\vert \right] \right\vert%
\end{array}%
\right.  \notag
\end{align}%
for any $x,y\in X.$

If we multiply (\ref{5.5}) by $p\left( x\right) p\left( y\right) \geq 0$ and
integrate, we deduce%
\begin{multline}
\int_{X}\int_{X}p\left( x\right) p\left( y\right) \left( \frac{q\left(
x\right) }{p\left( x\right) }-\frac{q\left( y\right) }{p\left( y\right) }%
\right) \left( g\left[ \frac{q\left( x\right) }{p\left( x\right) }\right] -g%
\left[ \frac{q\left( y\right) }{p\left( y\right) }\right] \right) d\mu
\left( x\right) d\mu \left( y\right)  \label{5.6} \\
\geq \left\{ 
\begin{array}{l}
\left\vert \int_{X}\int_{X}p\left( x\right) p\left( y\right) \left[
\left\vert \frac{q\left( x\right) }{p\left( x\right) }-1\right\vert
-\left\vert \frac{q\left( y\right) }{p\left( y\right) }-1\right\vert \right]
\right. \\ 
\ \ \ \ \ \ \ \times \left. \left[ g\left[ \frac{q\left( x\right) }{p\left(
x\right) }\right] -g\left[ \frac{q\left( y\right) }{p\left( y\right) }\right]
\right] d\mu \left( x\right) d\mu \left( y\right) \right\vert \\ 
\\ 
\left\vert \int_{X}\int_{X}p\left( x\right) p\left( y\right) \left[
\left\vert \frac{q\left( x\right) }{p\left( x\right) }-1\right\vert
-\left\vert \frac{q\left( y\right) }{p\left( y\right) }-1\right\vert \right]
\right. \\ 
\ \ \ \ \ \ \ \times \left. \left[ \left\vert g\left[ \frac{q\left( x\right) 
}{p\left( x\right) }\right] \right\vert -\left\vert g\left[ \frac{q\left(
y\right) }{p\left( y\right) }\right] \right\vert \right] d\mu \left(
x\right) d\mu \left( y\right) \right\vert%
\end{array}%
\right.
\end{multline}%
for any $x,y\in X.$

Now, observe that%
\begin{align*}
& \int_{X}\int_{X}p\left( x\right) p\left( y\right) \left[ \left\vert \frac{%
q\left( x\right) }{p\left( x\right) }-1\right\vert -\left\vert \frac{q\left(
y\right) }{p\left( y\right) }-1\right\vert \right] \left[ g\left[ \frac{%
q\left( x\right) }{p\left( x\right) }\right] -g\left[ \frac{q\left( y\right) 
}{p\left( y\right) }\right] \right] d\mu \left( x\right) d\mu \left( y\right)
\\
& =\int_{X}\int_{X}p\left( x\right) p\left( y\right) \left[ \left\vert \frac{%
q\left( x\right) }{p\left( x\right) }-1\right\vert g\left[ \frac{q\left(
x\right) }{p\left( x\right) }\right] +\left\vert \frac{q\left( y\right) }{%
p\left( y\right) }-1\right\vert g\left[ \frac{q\left( y\right) }{p\left(
y\right) }\right] \right] d\mu \left( x\right) d\mu \left( y\right) \\
& \quad -\int_{X}\int_{X}p\left( x\right) p\left( y\right) \left[ \left\vert 
\frac{q\left( x\right) }{p\left( x\right) }-1\right\vert g\left[ \frac{%
q\left( y\right) }{p\left( y\right) }\right] +\left\vert \frac{q\left(
y\right) }{p\left( y\right) }-1\right\vert g\left[ \frac{q\left( x\right) }{%
p\left( x\right) }\right] \right] d\mu \left( x\right) d\mu \left( y\right)
\\
& =2\int_{X}p\left( y\right) d\mu \left( y\right) \int_{X}p\left( x\right)
\left\vert \frac{q\left( x\right) }{p\left( x\right) }-1\right\vert g\left[ 
\frac{q\left( x\right) }{p\left( x\right) }\right] d\mu \left( x\right) \\
& \quad -2\int_{X}p\left( x\right) \left\vert \frac{q\left( x\right) }{%
p\left( x\right) }-1\right\vert d\mu \left( x\right) \int_{X}p\left(
y\right) g\left[ \frac{q\left( y\right) }{p\left( y\right) }\right] d\mu
\left( y\right) \\
& =2\left[ \bar{\delta}_{g}\left( Q,P\right) -V\left( Q,P\right) I_{g}\left(
Q,P\right) \right] ,
\end{align*}%
and a similar identity holds for the quantity in the second branch of (\ref%
{5.6}).

Finally, using the representation (\ref{5.1}), we deduce the desired
inequality (\ref{5.4}).
\end{proof}

\section{The Positivity of $\protect\delta -$Divergence for $g\in \mathcal{L}%
e_{1}\left( [0,\infty )\right) $}

The following result extending the positivity of $\delta -$divergence for
monotonic functions, holds.

\begin{theorem}
\label{t6.1}If $g\in \mathcal{L}e_{1}\left( [0,\infty )\right) ,$ then $%
\delta _{g}\left( Q,P\right) \geq 0$ for any $P,Q\in \mathcal{P}$.
\end{theorem}

\begin{proof}
We use the identity%
\begin{align}
\delta _{g}\left( Q,P\right) & =\int_{X}\left[ q\left( x\right) -p\left(
x\right) \right] g\left[ \frac{q\left( x\right) }{p\left( x\right) }\right]
d\mu \left( x\right)  \label{6.1} \\
& =\int_{X}p\left( x\right) \left[ \frac{q\left( x\right) }{p\left( x\right) 
}-1\right] g\left[ \frac{q\left( x\right) }{p\left( x\right) }\right] d\mu
\left( x\right)  \notag \\
& =\int_{X}p\left( x\right) \left[ \frac{q\left( x\right) }{p\left( x\right) 
}-1\right] \left[ g\left[ \frac{q\left( x\right) }{p\left( x\right) }\right]
-g\left( 1\right) \right] d\mu \left( x\right) .  \notag
\end{align}%
Since $g\in \mathcal{L}e_{1}\left( [0,\infty )\right) ,$ then for any $t\in
\lbrack 0,\infty )$ we have%
\begin{equation*}
\left( t-1\right) \left[ g\left( t\right) -g\left( 1\right) \right] \geq 0
\end{equation*}%
giving that%
\begin{equation*}
\left( \frac{q\left( x\right) }{p\left( x\right) }-1\right) \left[ g\left[ 
\frac{q\left( x\right) }{p\left( x\right) }\right] -g\left( 1\right) \right]
\geq 0
\end{equation*}%
for any $x\in X.$

Using the representation (\ref{6.1}), we deduce the desired result.
\end{proof}

\begin{corollary}
\label{c6.2}If $f\in \mathcal{O}_{0}\left( [0,\infty )\right) ,$ then $%
I_{f}\left( Q,P\right) \geq 0$ for any $P,Q\in \mathcal{P}$.
\end{corollary}

\begin{proof}
If $f\in \mathcal{O}_{0}\left( [0,\infty )\right) ,$ then there exists a $%
g\in \mathcal{L}e_{1}\left( [0,\infty )\right) $ such that $f\left( t\right)
=\left( t-1\right) g\left( t\right) $ for any $t\in \lbrack 0,\infty ).$ Then%
\begin{eqnarray*}
I_{f}\left( Q,P\right) &=&\int_{X}p\left( x\right) f\left[ \frac{q\left(
x\right) }{p\left( x\right) }\right] d\mu \left( x\right) \\
&=&\int_{X}p\left( x\right) \left[ \frac{q\left( x\right) }{p\left( x\right) 
}-1\right] g\left[ \frac{q\left( x\right) }{p\left( x\right) }\right] d\mu
\left( x\right) \\
&=&\delta _{g}\left( Q,P\right) \geq 0,
\end{eqnarray*}%
and the proof is completed.
\end{proof}

The following improvement of Theorem \ref{t6.1} holds.

\begin{theorem}
\label{t6.3}If $g\in \mathcal{L}e_{1}\left( [0,\infty )\right) ,$ then%
\begin{equation}
\delta _{g}\left( Q,P\right) \geq \left\vert \delta _{\left\vert
g\right\vert }\left( Q,P\right) \right\vert \geq 0  \label{6.2}
\end{equation}%
for any $P,Q\in \mathcal{P}$.
\end{theorem}

\begin{proof}
Since $g\in \mathcal{L}e_{1}\left( [0,\infty )\right) ,$ we obviously have%
\begin{align}
& \left[ \frac{q\left( x\right) }{p\left( x\right) }-1\right] \left[ g\left[ 
\frac{q\left( x\right) }{p\left( x\right) }\right] -g\left( 1\right) \right]
\label{6.3} \\
& =\left\vert \left( \frac{q\left( x\right) }{p\left( x\right) }-1\right)
\left( g\left[ \frac{q\left( x\right) }{p\left( x\right) }\right] -g\left(
1\right) \right) \right\vert  \notag \\
& \geq \left\vert \left( \frac{q\left( x\right) }{p\left( x\right) }%
-1\right) \left( \left\vert g\left[ \frac{q\left( x\right) }{p\left(
x\right) }\right] \right\vert -\left\vert g\left( 1\right) \right\vert
\right) \right\vert .  \notag
\end{align}%
Multiplying (\ref{6.3}) by $p\left( x\right) \geq 0$ and integrating on $X,$
we have%
\begin{align*}
& \int_{X}p\left( x\right) \left[ \frac{q\left( x\right) }{p\left( x\right) }%
-1\right] \left[ g\left[ \frac{q\left( x\right) }{p\left( x\right) }\right]
-g\left( 1\right) \right] d\mu \left( x\right) \\
& =\int_{X}p\left( x\right) \left\vert \left( \frac{q\left( x\right) }{%
p\left( x\right) }-1\right) \left( \left\vert g\left[ \frac{q\left( x\right) 
}{p\left( x\right) }\right] \right\vert -\left\vert g\left( 1\right)
\right\vert \right) \right\vert d\mu \left( x\right) \\
& \geq \left\vert \int_{X}p\left( x\right) \left( \frac{q\left( x\right) }{%
p\left( x\right) }-1\right) \left( \left\vert g\left[ \frac{q\left( x\right) 
}{p\left( x\right) }\right] \right\vert -\left\vert g\left( 1\right)
\right\vert \right) d\mu \left( x\right) \right\vert \\
& =\left\vert \delta _{\left\vert g\right\vert }\left( Q,P\right)
\right\vert ,
\end{align*}%
and the inequality (\ref{6.2}) is proved.
\end{proof}

\section{Bounds in Terms of the $\protect\chi ^{2}-$Divergence}

The following result may be stated.

\begin{theorem}
\label{t7.1}Let $g:\left[ 0,\infty \right] \rightarrow \mathbb{R}$ be a
differentiable function such that there exists the constants $\gamma ,\Gamma
\in \mathbb{R}$ with%
\begin{equation}
\gamma \leq g^{\prime }\left( t\right) \leq \Gamma \text{ \ \ for any \ }%
t\in \left( 0,\infty \right) .  \label{7.1}
\end{equation}%
Then we have the inequality%
\begin{equation}
\gamma D_{\chi ^{2}}\left( Q,P\right) \leq \delta _{g}\left( Q,P\right) \leq
\Gamma D_{\chi ^{2}}\left( Q,P\right) ,  \label{7.2}
\end{equation}%
for any $P,Q\in \mathcal{P}$.
\end{theorem}

\begin{proof}
Consider the auxiliary function $h_{\gamma }:\left[ 0,\infty \right]
\rightarrow \mathbb{R}$, $h_{\gamma }\left( t\right) :=g\left( t\right)
-\gamma \left( t-1\right) .$ Obviously, $h_{\gamma }$ is differentiable on $%
\left( 0,\infty \right) $ and since, by (\ref{7.1}),%
\begin{equation*}
h_{\gamma }^{\prime }\left( t\right) =g^{\prime }\left( t\right) -\gamma
\geq 0
\end{equation*}%
it follows that $h_{\gamma }$ is monotonic nondecreasing on $[0,\infty ).$

Applying Theorem \ref{t5.1}, we deduce%
\begin{equation*}
\delta _{h_{\gamma }}\left( Q,P\right) \geq 0\text{ \ \ for any \ }P,Q\in 
\mathcal{P}
\end{equation*}%
and since%
\begin{align*}
\delta _{h_{\gamma }}\left( Q,P\right) & =\delta _{g-\gamma \left( \cdot
-1\right) }\left( Q,P\right) \\
& =\int_{X}\left[ q\left( x\right) -p\left( x\right) \right] \left[ g\left[ 
\frac{q\left( x\right) }{p\left( x\right) }\right] -\gamma \left[ \frac{%
q\left( x\right) }{p\left( x\right) }-1\right] \right] d\mu \left( x\right)
\\
& =\delta _{g}\left( Q,P\right) -\gamma D_{\chi ^{2}}\left( Q,P\right) ,
\end{align*}%
then the first inequality in (\ref{7.2}) is proved.

The second inequality may be proven in a similar manner by using the
auxiliary function $h_{\Gamma }:\left[ 0,\infty \right) \rightarrow \mathbb{R%
}$, $h_{\Gamma }\left( t\right) :=\Gamma \left( t-1\right) -g\left( t\right)
.$
\end{proof}

The following corollary is a natural application of the above theorem.

\begin{corollary}
\label{c7.2}Let $f:\left[ 0,\infty \right] \rightarrow \mathbb{R}$ be a
differentiable convex function on $\left( 0,\infty \right) $ with $f\left(
1\right) =0.$ If there exist the constants $\gamma ,\Gamma \in \mathbb{R}$
with the property that:%
\begin{equation}
\gamma \left( t-1\right) ^{2}+f\left( t\right) \leq f^{\prime }\left(
t\right) \left( t-1\right) \leq f\left( t\right) +\Gamma \left( t-1\right)
^{2}  \label{7.3}
\end{equation}%
for any $t\in \left( 0,\infty \right) ,$ then we have the inequality:%
\begin{equation}
\gamma D_{\chi ^{2}}\left( Q,P\right) \leq I_{f}\left( Q,P\right) \leq
\Gamma D_{\chi ^{2}}\left( Q,P\right)  \label{7.4}
\end{equation}%
for any $P,Q\in \mathcal{P}$.
\end{corollary}

\begin{proof}
We know that for any $P,Q\in \mathcal{P}$, we have (see for example (\ref%
{4.1})):%
\begin{equation*}
I_{f}\left( Q,P\right) =\delta _{g_{f^{\prime }\left( 1\right) }}\left(
Q,P\right) ,
\end{equation*}%
where%
\begin{equation*}
g_{f^{\prime }\left( 1\right) }=\left\{ 
\begin{array}{ll}
\dfrac{f\left( t\right) }{t-1} & \text{if\ \ }t\in \lbrack 0,1)\cup \left(
1,\infty \right) , \\ 
&  \\ 
f^{\prime }\left( 1\right) & \text{if\ \ }t=1.%
\end{array}%
\right.
\end{equation*}%
We observe that, by the hypothesis of the corollary, $g_{f^{\prime }\left(
1\right) }$ is differentiable on $\left( 0,\infty \right) $ and 
\begin{equation*}
g_{f^{\prime }\left( 1\right) }^{\prime }\left( t\right) =\frac{f^{\prime
}\left( t\right) \left( t-1\right) -f\left( t\right) }{\left( t-1\right) ^{2}%
}
\end{equation*}%
for any $t\in \left( 0,1\right) \cup \left( 1,\infty \right) .$

Using (\ref{7.3}), we deduce that%
\begin{equation*}
\gamma \leq g_{f^{\prime }\left( 1\right) }^{\prime }\left( t\right) \leq
\Gamma
\end{equation*}%
for $t\in \left( 0,\infty \right) ,$ and applying Theorem \ref{t7.1} above,
for $g=g_{f^{\prime }\left( 1\right) },$ we deduce the desired inequality (%
\ref{7.4}).
\end{proof}

\section{Bounds in Terms of the $J-$Divergence}

We recall that the \textit{Jeffreys divergence} (or $J-$divergence for
short) is defined as%
\begin{equation}
J\left( Q,P\right) :=\int_{X}\left[ q\left( x\right) -p\left( x\right) %
\right] \ln \left[ \frac{q\left( x\right) }{p\left( x\right) }\right] d\mu
\left( x\right) ,  \label{8.1}
\end{equation}%
where $P,Q\in \mathcal{P}$.

The following result holds.

\begin{theorem}
\label{t8.1}Let $g:\left[ 0,\infty \right] \rightarrow \mathbb{R}$ be a
differentiable function such that there exists the constants $\phi ,\Phi \in 
\mathbb{R}$ with%
\begin{equation}
\phi \leq tg^{\prime }\left( t\right) \leq \Phi \text{ \ \ for any \ }t\in
\left( 0,\infty \right) .  \label{8.2}
\end{equation}%
Then we have the inequality%
\begin{equation}
\phi J\left( Q,P\right) \leq \delta _{g}\left( Q,P\right) \leq \Phi J\left(
Q,P\right) ,  \label{8.3}
\end{equation}%
for any $P,Q\in \mathcal{P}$.
\end{theorem}

\begin{proof}
Consider the auxiliary function $h_{\phi }:\left[ 0,\infty \right)
\rightarrow \mathbb{R}$, $h_{\phi }\left( t\right) :=g\left( t\right) -\phi
\ln t.$ Obviously, $h_{\phi }$ is differentiable on $\left( 0,\infty \right) 
$ and, by (\ref{8.2}),%
\begin{equation*}
h_{\phi }^{\prime }\left( t\right) =g^{\prime }\left( t\right) -\frac{\phi }{%
t}=\frac{1}{t}\left[ tg^{\prime }\left( t\right) -\phi \right] \geq 0,
\end{equation*}%
for any $t\in \left( 0,\infty \right) ,$ showing that the function is
monotonic nondecreasing on $(0,\infty ).$

Applying Theorem \ref{t5.1}, we deduce%
\begin{equation*}
\delta _{h_{\phi }}\left( Q,P\right) \geq 0\text{ \ \ for any \ }P,Q\in 
\mathcal{P}
\end{equation*}%
and since%
\begin{align*}
\delta _{h_{\phi }}\left( Q,P\right) & =\delta _{g-\phi \ln \left( \cdot
\right) }\left( Q,P\right) \\
& =\int_{X}\left[ q\left( x\right) -p\left( x\right) \right] \left[ g\left[ 
\frac{q\left( x\right) }{p\left( x\right) }\right] -\phi \ln \left[ \frac{%
q\left( x\right) }{p\left( x\right) }\right] \right] d\mu \left( x\right) \\
& =\delta _{g}\left( Q,P\right) -\phi J\left( Q,P\right) ,
\end{align*}%
then the first inequality in (\ref{8.3}) is proved.

The second inequality may be proven in a similar manner by using the
auxiliary function $h_{\Phi }:\left[ 0,\infty \right) \rightarrow \mathbb{R}$%
, $h_{\Phi }\left( t\right) :=\Phi \ln t-g\left( t\right) .$
\end{proof}

The following corollary is a natural application of the above theorem.

\begin{corollary}
\label{c8.2}Let $f:\left[ 0,\infty \right] \rightarrow \mathbb{R}$ be a
differentiable convex function on $\left( 0,\infty \right) $ with $f\left(
1\right) =0.$ If there exist the constants $\phi ,\Phi \in \mathbb{R}$ with
the property that:%
\begin{equation}
\phi \left( t-1\right) ^{2}+tf\left( t\right) \leq t\left( t-1\right)
f^{\prime }\left( t\right) \leq tf\left( t\right) +\Phi \left( t-1\right)
^{2}  \label{8.4}
\end{equation}%
for any $t\in \left( 0,\infty \right) ,$ then we have the inequality:%
\begin{equation}
\phi J\left( Q,P\right) \leq I_{f}\left( Q,P\right) \leq \Phi J\left(
Q,P\right)  \label{8.5}
\end{equation}%
for any $P,Q\in \mathcal{P}$.
\end{corollary}

The proof is similar to the one in Corollary \ref{c7.2} and we omit the
details.

\end{document}